\documentclass[12pt]{article}
\usepackage[applemac]{inputenc}
\usepackage[english]{babel}

\usepackage{amsmath} 
\usepackage{amssymb} 
\usepackage{ulem}
\usepackage{amsfonts} 
\usepackage{array} 

\def \R{\mathbb{R}}

\def\e{\varepsilon}

\def \Z{\mathbb{Z}}

\newtheorem{theorem}{Theorem}[section]

\newtheorem{Lemma}{Lemma}[section]

\newtheorem{Remark}{Remark}[section]

\def\lb{\overline\lambda}
\def\Fb{\overline F}
\def\wb{\overline w}
\def\yb{\overline y}
\def\Kb{\overline K}
\def\wa{w^\alpha}
\def\ou{\overline u}
\def\uu{\underline u}
\def\oU{\overline u}
\def\uU{\underline u}
\def\oC{\overline \chi}
\def\uC{\underline \chi}

\def\ue{u^{\varepsilon}}
\def\Ue{U^{\varepsilon}}
\def\Ce{\chi^{\varepsilon}}
\def\wa{w^\alpha}

\def\xb{\overline x}
\def\yb{\overline y}
\def\tb{\overline t}

\def\hyp#1{{\bf (H#1)}}
\def\1{1\!\!1}

\begin{document}
\title{Some Homogenization Results for Non-Coercive Hamilton-Jacobi Equations}
\author{Guy Barles
\\ \small Laboratoire de Math\'ematiques et Physique
Th\'eorique (UMR CNRS 6083) \\ \small
Fédération Denis Poisson \\ \small 
Universit\'e Fran\c{c}ois Rabelais de Tours \\ \small Parc de Grandmont, 37200 Tours, {\sc France} }
\date{ }
\maketitle

\begin{abstract} {\footnotesize Recently, C. Imbert \& R. Monneau study the homogenization of coercive Hamilton-Jacobi Equations with a $u/\e$-dependence : this unusual dependence leads to a non-standard cell problem and, in order to solve it, they introduce new ideas to obtain the estimates on the oscillations of the solutions. In this article, we use their ideas to provide new homogenization results for ``standard'' Hamilton-Jacobi Equations (i.e. without a $u/\e$-dependence) but in the case of {\it non-coercive Hamiltonians}. As a by-product, we obtain a simpler and more natural proof of the results of C. Imbert \& R. Monneau, but under slightly more restrictive assumptions on the Hamiltonians.} 
\end{abstract}

\footnote{
\noindent{\small \hspace{-0.6cm}\it{Key words and phrases :} {\rm Homogenization, Hamilton-Jacobi Equations, ergodic problems, level-set approach, viscosity solutions}}
}

\section*{Introduction}

In this article, we are interested in homogenization problems for first-order Hamilton-Jacobi Equations. The originality of this work is to provide results in the case of non-coercive Hamiltonians and applications to non-standard problems. Before describing more specifically our contributions, we want to point out that most of the new ideas used in this paper are borrowed from Imbert \& Monneau \cite{IM} who study the homogenization of (coercive) Hamilton-Jacobi Equations with a $\ue / \e$-dependence, namely
\begin{equation}\label{HJe}
\ue_t + H(\e^{-1}x, \e^{-1}\ue, D\ue)=0\quad \hbox{in}\  \R^n \times (0,+\infty),
\end{equation}
where $H$ is a continuous and coercive Hamiltonian. In fact, the starting point of the present work and one of its main motivations was to solve such problems, maybe under more restrictive assumptions than in \cite{IM} but (i) with simpler and more natural methods and (ii) with a clearer explanation of the involved phenomenas. We will come back later on (\ref{HJe}) and, in particular, we will explain the connections between the standard problems we are first dealing with and these non-standard problems.

In order to emphasize the main new ideas in our approach, we are not going to try to consider the most general framework but we restrict ourselves to a model case which carries the main difficulties. More precisely, we are going to study the limit as $\e \to 0$ of the solution $\Ue (x,y,t) $ of the Hamilton-Jacobi Equation
\begin{equation}\label{FJe}
\Ue_t + F(\e^{-1}x, \e^{-1}y, \e^{-1}t, D_x \Ue, D_y  \Ue)=0\quad \hbox{in}\  \R^{n+1} \times (0,+\infty),
\end{equation}
where $x \in \R^n$, $y\in \R$, $t \in (0,+\infty)$, $F(x,y,t,p_x,p_y)$ is a continuous function on $\R^{2n+3}$ which is $\Z^n$-periodic in $x$ and $1$-periodic in $y$ and $t$. More precise conditions on $F$ will be given later on but we point out that $F$ is assumed to be coercive with respect to $p_x$ but not with respect to $p_y$ and this is the key new point. In order to give a flavour of the type of Hamiltonians we are able to consider, a typical example is
$$ F(x,y,t,p_x,p_y)= a(x,t)|p_x|^\beta + b(x,y,t)|p_y| - f(x,t)\; ,$$
where $a,b, f$ are Lipschitz continuous function with the right periodicity in $x,y,t$ with $a(x,t) \geq \eta >0$ in $\R^n \times [0,+\infty)$ and $\beta \geq 1$. In particular, no sign condition is imposed on $b$, it may change sign and even the term $b(x,y,t)|p_y|$ can be replaced as well by a linear term $b(x,y,t)p_y$.

In order to understand the contribution of this work, we recall that the theory of homogenization for first-order Hamilton-Jacobi Equation started with the famous unpublished work of Lions, Papanicolaou \& Varadhan \cite{LPV} who completely solve the problem in the case of time-independent, periodic and \textit{coercive} Hamiltonians i.e. for the above equation when
$$ F(x,y, p_x,p_y) \to + \infty \hbox{  when }|p_x| + |p_y| \to + \infty,\hbox{ uniformly w.r.t $x$ and $y$.}$$
Then, to the best of our knowledge, such an assumption was used (one way or the other) in all the works concerning the homogenization of Hamilton-Jacobi Equations~: for more general periodic situations (cf. Ishii\cite{I1}), for problems set in bounded domains (cf. Alvarez\cite{A}, Horie \& Ishii\cite{HI}), for equation with different structure (cf. Alvarez \& Ishii\cite{AI}), for deterministic control problems in $L^\infty$ (cf. Alvarez \& Barron \cite{AB}), for almost periodic Hamiltonians (cf. Ishii\cite{I2}) and for Hamiltonians with stochastic dependence (cf. Souganidis\cite{PES}).

The reason why this coercivity assumption plays a central role, is that it allows to solve the so-called ``cell problem'' which provides the homogenized equation, i.e. the homogenized Hamiltonian. When $F$ is independent of time, this problem, which turns out to be an ergodic problem, consists in solving the pde
$$ F(x,y,D_x v + p_x,D_y v + p_y) = \lambda \quad\hbox{in  }\R^{n+1}\; ,$$
and in showing that, for any $(p_x,p_y)$, there exists a unique constant $\lambda=\Fb ( p_x , p_y)$ such that this equation has a bounded (periodic) solution. In general, solving this cell problem is the main difficulty and, in most cases, the fact that $\Fb$ is, indeed, the homogenized Hamiltonian, follows rather easily. It is worth pointing out, anyway, that, in order to have a comparison result for the limiting equation, one has to prove that $\Fb$ satisfies suitable properties and this may be a difficulty in some cases. We refer the reader to Section~\ref{ER} for a discussion in this direction.

In this framework, the role of the coercivity is rather clear since it provides an easy bound on $Dv = (D_x v, D_y v)$ once one knows that $ F(x,y,D_x v + p_x,D_y v + p_y)$ is bounded and this is a key argument to solve the cell problem. Up to now (again, to the best of our knowledge), no work (except perhaps partially \cite{PES}) succeeded to find an argument to solve this cell problem by pde methods in a general way which avoids the coercivity assumption on $F$ (and actually one of the main question is to find the right set of assumptions to do it). Furthermore, when $F$ depends on $t$, one has to find a space-time periodic solution of 
\begin{equation}\label{CP}
v_t + F(x,y,t,D_x v + p_x,D_y v + p_y) = \Fb ( p_x , p_y) \quad\hbox{in  }\R^{n+1} \times \R\; ,
\end{equation}
and the assumptions to do it seem even stronger (cf. Souganidis and the author \cite{BS}).

However it is worth pointing out that the pure pde approach we describe above for both homogenization and ergodic problems can be replaced, in some cases, by an approach using either optimal control (if $F$ is convex) or differential games methods : in this framework, results do exist for non-coercive Hamilton-Jacobi Equations, the coercivity assumption being replaced by either controllability or non-resonance assumptions. We refer the reader to Alvarez and Bardi \cite{AB1, AB2, AB3}, Arisawa \cite{A1,A2}, Arisawa and Lions \cite{AL} and Artstein and Gaitsgory \cite{AG}, Bardi \cite{B} for results in this direction. 

In order to solve problems like (\ref{CP}) in a general way and once for all $(p_x,p_y)$, we consider the ergodic problem~: find a constant $\lb$ such that the equation
\begin{equation}\label{EP}
 w_t + G (x,y,t,D_x w ,D_y w) = \lb \quad\hbox{in  }\R^{n+1} \times \R\; ,
\end{equation}
has a bounded (periodic) solution. In fact, with our approach, we obtain a weaker but sufficient result which can take two possible forms~: we can prove that there exists a unique constant $\lb$ such that we have either \textit{approximate} continuous (but not necessarily Lipschitz continuous) periodic solutions or \textit{exact}, possibly discontinuous, periodic sub and supersolutions. This result is the first main contribution of the paper and this is where we use in a key way the ideas of Imbert \& Monneau \cite{IM} ; an unusual feature of the proof is the estimate of the oscillation $\max_{\R^{n+1} \times \R} (w) - \min_{\R^{n+1} \times \R} (w) $ which replaces the classical gradient estimate and relies on two very original arguments.

Once (\ref{EP}) is solved, the result for the homogenization problem follows by using the usual arguments~: indeed, it is enough to (essentially) apply the result to (\ref{CP}) for any $(p_x,p_y)$, then to obtain suitable properties for $\Fb$ and finally to use the ``perturbed test function's method'' of Evans \cite{LCE1,LCE2}, even if, to prove the convergence, we have to introduce an additional argument to take care of the rather weak properties we impose on $F$ in the variables $x$, $t$ and $p_x$. We refer the reader to the papers of Alvarez \& Bardi \cite{AB1,AB2} for a general approach of singular perturbation problems for fully nonlinear, possibly degenerate, elliptic or parabolic equations; this approach applies to homogenization problems, including in the case of first-order Hamilton-Jacobi Equations, and, in particular, the authors clarify the connections between ergodic and homogenization problems in a general setting.

Now we come back to (\ref{HJe}) which seems to have nothing to do with the above (rather) classical problem. Roughly speaking, the connection consists in considering (\ref{HJe}) as the equation of the motion of the graph $y=\ue (x,t)$ and to introduce the associated level set equation which is an equation in $\R^{n+1}$. Since, intuitively, the level set equation should play a role only on the moving hypersurface, this leads to replace $\e^{-1}\ue(x,t)$ by $\e^{-1}y$ and to transform (\ref{HJe}) in a standard (but a priori non-coercive) homogenization problem. The rigourous justification of this formal (but convincing) argument is easy by following the article of Giga \& Sato \cite{GS} which is exactly describing the way to do it for first-order equations (see also the paper of Biton, Ley \& the author\cite{BBL} for motion by Mean Curvature). Compared to the article of Imbert \& Monneau \cite{IM}, this approach leads to more restrictive assumptions on $H$ but to a far simpler proof and to more natural and understandable arguments. Typically the results of \cite{IM} apply for all $H$ of the form
$$ H(x,y,t,p) = a(x,t) (|p|^2+1)^{\beta/2} + g(x,y,t) \; ,$$
where $a, g$ are continuous periodic functions with $a>0$ in $\R^n \times \R$ and $0< \beta \leq 1$, while here we have to consider the case $\beta = 1$. Finally, it is worth pointing out that Imbert \& Monneau \cite{IM} use also an extension to $\R^{n+1}$ by adding a extra variable but with no clear interpretation of this new variable which is just a trick in the proof.

The paper is organized as follows~: in the first section, we study the ergodic problem (\ref{EP}); this is the occasion to present the key assumptions on $G$ (and therefore on $F$) which describe the structure of the problem allowing the lack of coercivity. In Section~\ref{Hom}, we apply this result to the cell problem and we solve the homogenization problem; to do so, we show that $\Fb ( p_x , p_y)$ is a continuous function, which depends also continuously on $F$. Then, in Section~\ref{Homue}, we consider the problem (\ref{HJe}): we describe the Giga \& Sato \cite{GS} approach and the assumptions on $H$ to apply it; the conclusion then follows from the results of Section~\ref{Hom}. Finally, we provide remarks on the possible variants and extensions of the above results in Section~\ref{ER}.

\medskip
\noindent{\bf Acknowledgments : }The author wishes to thank Martino Bardi for interesting comments on the first version of this paper and, in particular, for pointing out some missing references. He also thanks the anonymous referees for their careful reading of the paper and their very constructive remarks.

\section{The non-coercive ergodic problem}\label{sec:EP}

In order to state the result, we have to impose two kinds of assumptions: first some basic assumptions which ensure existence, uniqueness and comparison properties for the problems we are going to introduce and then structure conditions which allow to obtain the estimates we need to solve the ergodic problem. Again, to focus on the main ideas of the proof, we are not going to provide the most sophisticated assumptions but the simplest relevant ones.

\medskip
\noindent \hyp{1} $G(x,y,t,p_x,p_y)$ is a continuous function in $\R^{2n+3}$, $\Z^n$-periodic in $x$, $1$-periodic in $y$ and $t$.

\medskip
\noindent\hyp{2} Either $G$ is independent of $t$ and $ G (x,y,p_x,p_y) \to + \infty $ when $|p_x|\to +\infty$, uniformly w.r.t. $(x,y) \in \R^n \times \R $ and $|p_y| \leq R$ for any $R>0$, or $G$ depends on $t$ and there exists constants $C_1,C_2 >0$ such that
$$  G(x,y,t,p_x,0) \geq C_1 |p_x| - C_2 \quad\hbox{in  }\R^n \times \R \times \R \times \R^n.$$

\medskip
\noindent\hyp{3} The function $G$ is locally Lipschitz continuous w.r.t. $y$, $t$, $p_y$ and there exists $l \in \R$ and $C_3, C_4, C_5\geq 0$ such that, for almost every $\xi = (x,y,t,p_x,p_y) \in \R^{2n+3}$
$$  |D_y G (\xi) | \leq C_3 |p_y + l |\; ,\; |D_t G (\xi) | \leq C_4 (1+ |p_y| + |G (\xi)|)\; , \;  |D_{p_y} G (\xi) | \leq C_5 .$$

\medskip
Our result is the following
\begin{theorem}\label{thm:EP} Under assumptions \hyp{1}-\hyp{3}, we have\\
(i) For any $\alpha >0$, there exists a unique continuous, space-time periodic solution $\wa $ of
$$  \wa_t + G (x,y,t,D_x \wa ,D_y \wa) + \alpha \wa= 0 \quad\hbox{in  }\R^{n+1} \times \R\; ,$$
which is independent of $t$ if $G$ is independent of $t$.\\
(ii) There exists a constant $K= K(G)$ depending only on $G$ through
$C_0 := ||G(x,y,t,0,0)||_\infty$, $l$ and the coercivity assumption \hyp{2}, such that
$$ \max_{\R^{n+1} \times \R}\,\wa -  \min_{\R^{n+1} \times \R}\,\wa \leq K\; .$$
In particular, if $G$ depends on $t$, then $K= K(C_0, C_1,C_2,l)$.\\
(iii) There exists a unique constant $\lb$ such that the equation
\begin{equation}\label{ergeqn}
 w_t + G (x,y,t,D_x w ,D_y w) = \lb \quad\hbox{in  }\R^{n+1} \times \R\; ,
\end{equation}
admits bounded, possibly discontinuous, sub and supersolutions.\\
(iv) If $v_0$ is bounded uniformly continuous in $\R^{n+1}$, there exists a unique solution $v$ of the initial value problem
\begin{equation}\label{evol}
 v_t + G (x,y,t,D_x v ,D_y v) = 0 \quad\hbox{in  }\R^{n+1} \times (0,+\infty)\; ,
\end{equation}
\begin{equation}\label{evolid}
v(x,y,0) = v_0 (x,y) \quad\hbox{in  }\R^{n+1}\; ,
\end{equation}
which is bounded and uniformly continuous in $\R^{n+1} \times [0,T]$ for all $T>0$; moreover we have
$$ \lim_{t\to +\infty}\,t^{-1} v(x,y,t) = - \lb\; ,$$
the limit being uniform in $\R^{n+1}$.
\end{theorem}

\medskip
We have chosen to state Theorem~\ref{thm:EP} in that way for several reasons : first, even if the result (i) seems to be classical, it is not the case because the assumptions \hyp{1}-\hyp{3} are not completely standard and it may not be so clear that the $\wa$-equation has a unique solution. Indeed, we only assume $G$ to be continuous in $x$, which is a very weak assumption, but this is compensated by the coercivity of $G$ in $p_x$. To prove (i) which is the first (slight) difficulty of the proof, we borrow arguments from P.L. Lions and the author \cite{GBPLL} and the author \cite{BG1}. Then (ii) is the main step but (curiously) not the most technical one. (iii) is one of the possible conclusion but if one insists to use continuous ``correctors'', then (ii) implies that, for any $\delta >0$, if $\alpha$ is small enough, we have
$$ \lb - \delta \leq \wa_t + G (x,y,t,D_x \wa ,D_y \wa) \leq \lb + \delta \quad\hbox{in  }\R^{n+1} \times \R\; ,$$
since, as usual and by using (ii), $\lb$ is the uniform limit of $-\alpha \wa$. Finally (iv) is an other classical characterization of the ergodic constant $\lb$.

\medskip
\noindent{\bf Proof : }We provide the proof in the case when $G$ actually depends on $t$, the other case being simpler (see the remark at the end of the section). 

\noindent{\bf 1.} We start by (i). We just sketch the proof and refer to \cite{GBPLL}, \cite{BG1} and \cite{GBLivre} for the (easy) details. As it is classical in viscosity solutions' theory, the existence of $\wa$ relies on having a Strong Comparison Result for the equation, the existence being obtained by the Perron's method extended to the framework of viscosity solutions by Ishii \cite{I3}. Here, because of the periodicity of $G$, Perron's method can be applied in the set of periodic subsolutions (by remarking that if $w$ is as a subsolution, $ \bar w (x,y,t) := \sup\{ w(x+l,y+k,t+m,\; l \in \Z^n, k,m \in \Z\}$ is still a subsolution), we can prove the comparison result only for periodic sub and supersolution which avoids problems with the unboundedness of the domain.

To prove this comparison result, we argue in the following way: if $u$ is an usc periodic subsolution and $v$ a lsc periodic supersolution of the equation, we approximate $\max_{\R^{n+2}}\, (u-v) $ by $\max_{\R^{2n+4}}\, (\Psi) $ where $\Psi=\Psi_{\eta, \beta}$ is given by
$$ \Psi (x_1,x_2,y_1,y_2,t,s): = u(x_1,y_1,t) - v(x_2,y_2,s) - \frac{(x_1-x_2)^2}{\eta^2} -  \frac{(y_1-y_2)^2}{\beta^2} - \frac{(t-s)^2}{\beta^2}\; ,$$
where $\eta, \beta >0$ are small constants devoted to tend to $0$.

Because of the periodicity of $u$ and $v$, the usc function $\Psi$ achieves its maximum at some point, which is still denoted $(x_1,x_2,y_1,y_2,t,s)$ for the sake of simplicity of notations (but which depends on $\eta$ and $\beta$). We set
$$  p_x := \frac{2(x_1-x_2)}{\eta^2} \; , \;  p_y := \frac{2(y_1-y_2)}{\beta^2} \; , \; p_t := \frac{2(t-s)}{\beta^2}\; .$$

The viscosity solutions inequalities for $u$ and $v$ read
\begin{equation}\label{inequ}
p_t + G (x_1,y_1,t,p_x ,p_y ) + \alpha u(x_1,y_1,t) \leq 0\; ,
\end{equation}
$$ p_t + G (x_2,y_2,s,p_x ,p_y ) + \alpha v(x_2,y_2,s) \geq 0\; .$$
At this point we recall the well-known fact that 
 \begin{equation}\label{penter}
\frac{(x_1-x_2)^2}{\eta^2} +  \frac{(y_1-y_2)^2}{\beta^2} + \frac{(t-s)^2}{\beta^2} \to 0\quad\hbox{as  }\beta, \eta \to 0\; .
\end{equation}
(We point out here that this is true since we deal with time-space periodic sub and supersolutions). 

This implies, in particular, that $p_t = o(\beta^{-1})$ and since $u$ and $v$ are bounded, we may as well assume that, for $\beta$ small enough, the above viscosity inequalities holds with $G^\beta = \min(\beta^{-1}, \max (G,-\beta^{-1}))$, which has the same properties as $G$ except that $|D_t G^\beta (\xi)| \leq   C_4 (1+|p_y| + \beta^{-1}). $

Subtracting the above viscosity sub and supersolutions inequalities and using \hyp{3} then yield
\begin{eqnarray*}
 \alpha (u(x_1,y_1,t) - v(x_2,y_2,s)) & \leq &  G^\beta (x_2,y_1,t,p_x ,p_y ) - G^\beta (x_1,y_1,t,p_x ,p_y ) \\
 & & + C_3|p_y + l||y_1-y_2| + C_4|t-s| (1+|p_y| + \beta^{-1})\; .
\end{eqnarray*}

The argument is then the following : we first fix $\beta$ and let $\eta$ tend to $0$. Because of \hyp{2} and inequality (\ref{inequ}), the $x$-gradient term $p_x$ remains bounded since $p_t, p_y$ remains bounded~; furthermore, since $u$ and $v$ are periodic, we may also assume that $(x_1,x_2,y_1,y_2,t,s)$ remains in a compact subset of $\R^{2n+4}$; therefore the passage to the limit $\eta \to 0$ can be done by using only the continuity of $G^\beta$ (or, equivalently, of $G$), all terms being convergent up to a subsequence.

Therefore, we are left with
$$ \alpha (u(x_1,y_1,t) - v(x_2,y_2,s)) \leq C_3|p_y + l||y_1-y_2| + C_4|t-s| (1+|p_y| + \beta^{-1})\; , $$
and the conclusion follows by letting $\beta$ tend to $0$ : the terms of the right-hand side converge to $0$ because of (\ref{penter}), while we know that $u(x_1,y_1,t) - v(x_2,y_2,s) \to \max_{\R^{n+2}}\, (u-v)$; we reach the conclusion that $\max_{\R^{n+2}}\, (u-v) \leq 0$ and the proof of (i) is complete.

\noindent{\bf 2.} We further remark that, either by construction or by comparison with constants, we have 
\begin{equation}\label{estlb}
\min_{x,y,t}\, [-G (x,y,t,0,0)] \leq \alpha \wa \leq \max_{x,y,t}\, [-G(x,y,t,0,0)]\; .
\end{equation}
In particular, $||\alpha \wa||_\infty \leq C_0$.

\noindent{\bf 3.} Now we turn to (ii). The idea of Imbert \& Monneau consists in estimating separately the behavior of $\wa$ in $(x,t)$ and in $y$. For the estimate in $(x,t)$, we introduce the function
$$ \wb (x,t):=\max_{y\in \R}\, \wa(x,y,t)\; ,$$
which, using \hyp{2} and the estimate on $\alpha \wa$ given in Step~2, is a space-time periodic subsolution of
$$ \wb_t + C_1| D_x \wb| - (C_0+C_2) \leq 0 \quad\hbox{in  }\R^n \times \R\; .$$
But the Oleinik-Lax Formula shows that, for any $x,t$ and for any $s\leq t$
$$ \wb (x,t) \leq \min_{|z-x| \leq C_1(t-s)}\,[ \wb(z,s)] + (C_0+C_2) (t-s)\; ,$$ and therefore, for any $S>0$
$$ \wb (x,t) \leq \min_{0\leq t- s\leq S} \left\{\min_{|z-x| \leq C_1(t-s)}\,[ \wb(z,s)] + (C_0+C_2) (t-s)\right\}\; .$$
Choosing $(x,t)$ such that $ \wb (x,t) = \max_{\R^n \times \R}\,\wb$ and choosing $S$ large enough in order that the set $\{(z,s); 0\leq t- s\leq S\ ,\ |z-x| \leq C_1(t-s)\}$ contains a whole period, we deduce that
\begin{equation}\label{k1b}
 \max_{\R^n \times \R}\,\wb - \min_{\R^n \times \R}\, \wb \leq \Kb_1\; ,
\end{equation}
where $\Kb_1$ depends only on $C_0$, $C_1$ and $C_2$.

\noindent{\bf 4.} Next we consider the behavior of $\wa$ in $y$. We are going to show that $D_y (\wa (x,y,t)+ly)$ has the same sign as $l$. We only do it in the case $l<0$, the case $l>0$ being treated in an analogous way. If $l=0$ then $\wa$ is independent of $y$~: indeed, since $G(x,y,t,p_x,0)$ is independent of $y$ by \hyp{3}, by using the same arguments as in the proof of (i), one can build a solution for this new Hamiltonian which depends only on $x$ and $t$ and which is also a solution of the $G$ equation. Therefore, by uniqueness, $\wa=\wa(x,t)$.

To prove that $y \mapsto \wa (x,y,t)+ly$ is nonincreasing for $l<0$, we consider
$$ M:=\max\{ \wa (x,y_1,t) - \wa (x,y_2,t) + l(y_1 -y_2);\ x \in \R^n,\ t \in \R,\ y_1 \geq y_2\}\; .$$
We have to show that this maximum is nonpositive. To do so, we argue by contradiction assuming that $M>0$ and we introduce the function\\
$\displaystyle \chi (x_1,x_2,y_1,y_2,t,s): = \wa (x_1,y_1,t) - \wa (x_2,y_2,s) - \frac{(x_1-x_2)^2}{\eta^2} -  \frac{[(y_1-y_2)^-]^2}{\beta^2}$

{\hfill $\displaystyle  + l(y_1-y_2) - \frac{(t-s)^2}{\beta^2}\; .$}

The function $\chi$ has indeed a maximum since $\wa$ is continuous and periodic (hence bounded) and since $l$ is strictly negative, the term  $l(y_1-y_2)$ controls $(y_1-y_2)^+$ which therefore remains bounded.

Denoting by $(x_1,x_2,y_1,y_2,t,s)$ a maximum point of $\chi$ and
$$ a= \frac{2(t-s)}{\beta^2} \; , \;    p_x = \frac{2(x_1-x_2)}{\eta^2}   \; , \; p_y = -l +\frac{2(y_1-y_2)^-}{\beta^2} \; ,$$
we have the viscosity inequalities
$$ a + G(x_1,y_1,t,p_x,p_y) + \alpha \wa (x_1,y_1,t) \leq 0\; ,$$
$$ a + G(x_2,y_2,s,p_x,p_y) + \alpha \wa (x_2,y_2,s) \geq 0\; .$$
Using \hyp{3}, we have
$$ | G(x_1,y_1,t,p_x,p_y) - G(x_1,y_2,t,p_x,p_y)| \leq 2C_3  \frac{[(y_1-y_2)^-]^2}{\beta^2} \; ,$$
and subtracting the above inequalities, we are left with\\
$ \alpha (\wa (x_1,y_1,t) - \wa (x_2,y_2,s) + l(y_1 -y_2) ) \leq $

{\hfill $\displaystyle G(x_2,y_2,s,p_x,p_y)-G(x_1,y_2,t,p_x,p_y) + 2C_3  \frac{[(y_1-y_2)^-]^2}{\beta^2} + \alpha l(y_1 -y_2)\; .$}

To conclude, we argue as in the proof of the comparison result : first, by periodicity we may assume that $(x_1,y_1,t)$ and $(x_2,y_2,s)$ remain in a compact subset of $\R^n \times \R \times \R$; next, by classical arguments, all the penalization terms tend to $0$. In particular, $(y_1-y_2)^- = o(\beta)$ and this yields $ \alpha l(y_1 -y_2) \leq o(\beta)$ since $l<0$.
Therefore we first fix $\beta$ and let $\eta$ tend to $0$ and then we let $\beta$ tend to $0$, using the same arguments as above. This gives the inequality $ \alpha M \leq 0$ and the conclusion.

\noindent{\bf 5.} To prove (ii), it is enough to put together the two above informations and again we do it only for $l<0$: we consider two points $(x,y,t)$ and $(x',y',t')$ such that $\wa (x,y,t) = \max_{\R^{n+2}} \,\wa$ and $\wa (x',y',t') = \min_{\R^{n+2}} \, \wa$, and we denote by $\yb$ the smallest point $z\geq y'$ such that $\wb (x',t') = \wa (x', z, t')$; we recall that $\wa$ is periodic and therefore $\yb - y' \leq 1$. We have
\begin{eqnarray*}
\wa (x,y,t) & = & \wb (x,t) \\
& \leq & \wb (x',t') + \Kb_1 \\
& \leq & \wa (x', \yb, t') + \Kb_1 \\
& \leq &  \wa (x', y',t' ) + l (y' - \yb) + \Kb_1 \\
& \leq &  \wa (x', y',t' ) - l + \Kb_1\; .
\end{eqnarray*}
which concludes the proof of (ii).

\noindent{\bf 6.} For the proof of (iii), we notice that $\wa - \min\,\wa$ is bounded and so is $\alpha \min\,\wa$ since $||\alpha \wa||_\infty \leq C_0$. Without loss of generality, we may assume that $-\alpha \min\,\wa$ converges to a constant $\lb$ and, along the same subsequence, the half-relaxed limits of $\wa - \min\,\wa$, denoted by $\overline w$ and $\underline w$, provide respectively a bounded subsolution and a bounded supersolution of (\ref{ergeqn}).

Now we show that there exists a unique constant $\lb$ such that (\ref{ergeqn}) has a bounded subsolution and a bounded supersolution. First, we assume that there exists $\tilde \lambda$ such that (\ref{ergeqn}) has a bounded subsolution $\tilde u$. By comparison result for the initial value problem (a result proved in a similar way as in the stationary case), using that $\tilde u - \tilde \lambda t$ and $\underline w - \lb t$ are respectively sub and supersolution of equation (\ref{evol}), we deduce
$$ \max_{\R^{n+2}}\,(\tilde u - \underline w)(x,y,t) \leq \max_{\R^{n+2}}\,(\tilde u - \underline w)(x,y,0) + t(\tilde \lambda - \lb)\; .$$
Since $\tilde u$ and $\underline w$ are bounded, by dividing by $t>0$ and letting $t$ tend to $+\infty$, we obtain $\tilde \lambda - \lb \geq 0$. The reverse inequality is obtained by the same argument assuming that there exists a bounded supersolution for the same constant.

\noindent{\bf 7.} Finally we prove (iv). The proof of the existence and uniqueness of the solution $v$ follows along the line of the proof of (i) with few additional classical arguments and therefore we skip it. In order to prove the convergence property, we use $\overline w$ and $\underline w$ defined in step 6 above. Since $v_0$,  $\overline w$ and $\underline w$ are bounded, there exist constants $\underline c, \overline c$ such that
$$ \underline c + \overline w (x,y,0) \leq v_0(x,y) \leq  \overline c + \underline w(x,y,0) \quad \hbox{in  }\R^{n+1}\; .$$
But $ \underline c + \overline w (x,y,t)- \lb t$ and $\overline c + \underline w(x,y,t) -\lb t$ are respectively sub and supersolution of (\ref{evol}) and, by comparison, we deduce
$$ \underline c + \overline w (x,y,t)- \lb t \leq v(x,y,t) \leq \overline c + \underline w(x,y,t) -\lb t \quad \hbox{in  }\R^{n+1}\times (0,+\infty)\; .$$
Dividing by $t$ and using the fact that $\overline w$ and $\underline w$ are bounded provide the property we wanted to prove.
And the proof of Theorem~\ref{thm:EP} is complete.

\begin{Remark} In the case where $G$ is independent of $t$, $\wa$ is also independent of $t$ since $\wa (x,y,t)$ and $\wa (x,y,t+h)$ are both solutions of the equation for any $h \in \R$ and therefore, by uniqueness, $\wa (x,y,t)=\wa (x,y,t+h)$ for any $h \in \R$. But maybe a more simple approach in this case is to solve the stationary equation directly and not to consider the evolution equation. In the above proof, the only change concerns the estimate obtained through $\wb=\wb (x)$ : indeed $\wb$ is a subsolution of
$$ \min_{y\in \R}\, G(x,y, D\wb, 0) \leq C_0 \quad \hbox{in  }\R^n\; ,$$
and this inequality together with the coercivity asssumption on $G$, shows that $||D\wb||_\infty \leq K$, where $K$ depends only on $G$ through $C_0$ and the coercivity. The function $\wb$ being $\Z^n$-periodic in $x$, this gives the bound on the oscillation and even a more precise information.
\end{Remark}

\section{Homogenization of the non-coercive equation}\label{Hom}

In order to state and prove the homogenization result, we first need to identify the homogenized Hamiltonian and its properties.

\begin{Lemma}\label{lem:Fb} Assume that $F$ satisfies the assumptions \hyp{1}-\hyp{3} with $l=0$ in \hyp{3}. For any $(p_x,p_y) \in \R^n \times \R$, there exists a unique constant $\Fb ( p_x , p_y) $ such that (\ref{CP}) has a bounded space-time periodic solution. The function $\Fb ( p_x , p_y) $ is a continuous function of $(p_x,p_y)$. Moreover, if $(F_k)_k$ is a sequence of functions satisfying the same assumptions as $F$ and uniformly in $k$, and which converges locally uniformly to $F$ on $\R^{2n +3}$, then $\Fb^k \to \Fb$ locally uniformly in $\R^n \times \R$.
\end{Lemma}

\medskip
\noindent{\bf Proof : }The first part of the result is an immediate consequence of Theorem~\ref{thm:EP} applied to
$$ G(x,y,t,q_x,q_y) = F(x,y,t,q_x + p_x,q_y+p_y)\; .$$

The continuity of $\Fb$ comes from all the estimates we have~: first we have a bound on $\Fb ( p_x , p_y)$ which is coming from (\ref{estlb}) (recall that $-\alpha \wa$ converges uniformly to $\lb$), namely
$$  \min_{x,y,t}\, [F(x,y,t, p_x,p_y)] \leq \Fb ( p_x , p_y) \leq \max_{x,y,t}\, [F(x,y,t, p_x,p_y)]\; .$$
Then we have also the bounds on the oscillations of the associated sub and supersolution $\overline v$, $\underline v$, typically through $\Kb_1$ given by (\ref{k1b}) in the time dependent case, which depend only on $( p_x , p_y)$, $F$ and $l=p_y$. 

Indeed, if the sequence $((p^k_x,p^k_y))_k$ is converging to $(p_x,p_y)$, then, up to a subsequence, $\Fb (p^k_x,p^k_y)$ converges to $\mu \in \R$. But if $\overline v^k$ and $\underline v^k$ are bounded sub and supersolution associated respectively to $(p^k_x,p^k_y)$, then they are uniformly bounded and the half-relaxed limit method provides $\overline v$ and $\underline v$ which are the sub and supersolution associated to $(p_x,p_y)$ and the constant $\mu$. By the uniqueness result of Theorem~\ref{thm:EP} (iii), this implies $\mu = \Fb ( p_x , p_y)$ and since this is true for any converging subsequence of $(\Fb (p^k_x,p^k_y))_k$, the continuity is proved. 

Finally the proof for the $\Fb^k$ relies on the same type of arguments and therefore we skip it.

\medskip
The homogenization result is then the
\begin{theorem} Assume that $F$ satisfies the assumptions \hyp{1}-\hyp{3} with $l=0$ in \hyp{3}. If the $\Ue$ are continuous viscosity solutions of (\ref{FJe}) associated with the initial data $\Ue (x,y,0) = U_0 (x,y)$ in $\R^n \times \R$ where $U_0$ is bounded uniformly continuous in $\R^{n+1}$, then, as $\e \to 0$, $\Ue \to U$ locally uniformly in $\R^{n+1} \times (0,+\infty)$, where $U$ is the unique solution of
\begin{equation}\label{homeqn}
U_t + \Fb (D_x U, D_y U) = 0 \quad \hbox{in}\  \R^{n+1} \times (0,+\infty),
\end{equation}
with the same initial datum $U_0$
\end{theorem}

\medskip
\noindent{\bf Proof :} We just sketch it since it follows essentially by classical arguments. We first remark that, if $A:=\max (|| U_0||_\infty , ||F(x,y,t,0,0) ||_\infty)$, then $-A(1+t)$ and $A(1+t)$ are respectively sub and supersolution of the problem and, by comparison, we have
$$ -A(1+t) \leq \Ue (x,y,t) \leq A(1+t) \quad \hbox{in  }\R^{n+1}\times (0,+\infty)\; ;$$
therefore,  the $\Ue$'s are uniformly bounded. 
  
Then we use the half-relaxed limit method and set $\oU := {\limsup}^* \,\Ue$, $\uU=\liminf_*\, \Ue$. We show how to prove that $\oU$ is a subsolution of (\ref{homeqn}), the proof that $\uU$ is a supersolution being analogous.

If $\phi$ is a smooth function and $(\xb, \yb, \tb)$ is a strict maximum point of $\oU - \phi$, we use the perturbed test-function's method in the following way: we set $p_x = D_x \phi (\xb, \yb, \tb)$, $p_y = D_y \phi (\xb, \yb, \tb)$ and let $V$ be a bounded space-time periodic supersolution of (\ref{CP}) associated to $(p_x,p_y)$ and to the function
$$ F^k (x,y,t,q_x,q_y) = \min\,\{F (x',y,t',q'_x,q_y)\ ;\ |x'-x| , |t'-t| , |q'_x-q_x|\leq k^{-1}\} \,, $$
for $k\geq 1$.

We introduce the function
$$ \Ue (x,y,t) - \phi (x,y,t) - \e V (\e^{-1}x', \e^{-1}y',\e^{-1}t') - \frac{|x-x'|^2}{\delta^2} - \frac{|y-y'|^2}{\delta^2}- \frac{|t-t'|^2}{\delta^2}\; .$$
It is easy to see (and standard) that, as $\e \to 0$ and $\delta \to 0$ with $\delta \ll \e$, the maximum points of this function converge to $(\xb, \yb, \tb,\xb, \yb, \tb)$. Performing standard computations (which we recall in the Appendix for the reader's convenience) with fixed $k$ leads to the inequality
\begin{equation}\label{fkvi} \phi_t (\xb, \yb, \tb) +  \Fb^k (D_x \phi (\xb, \yb, \tb),D_y \phi (\xb, \yb, \tb)) \leq 0\; ,
\end{equation}
and we use Lemma~\ref{lem:Fb} to conclude by letting $k\to\infty$. 

We just point out that the $F^k$-trick allows to take into account the rather weak conditions we impose on $ F (x,y,t,q_x,q_y) $ in $x, t$ and $q_x$~: indeed, if $(x,y,t,x',y',t')$ is a maximum point of the above function, we have, on one hand, $|D_x \phi (x,y,t)-D_x \phi (\xb, \yb, \tb)| \leq k^{-1}$ for $\e$ and $\delta$ small enough because $(x,y,t) \to (\xb,\yb,\tb)$ when $\e, \delta \to 0$; on the other hand, since the terms $\frac{|x-x'|^2}{\delta^2}$, $\frac{|t-t'|^2}{\delta^2}$ are bounded, then $|\e^{-1}x'-\e^{-1}x|, |\e^{-1}t' - \e^{-1}t| \leq k^{-1}$ if $\delta \ll \e$ is small enough. Therefore
$$ F^k (\e^{-1}x',\e^{-1}y',\e^{-1}t', D_x \phi (\xb, \yb, \tb) + q_x ,q_y) \leq F (\e^{-1}x,\e^{-1}y',\e^{-1}t,D_x \phi (x,y,t) + q_x ,q_y)\; ,$$ which is the inequality needed in the proof because the dependence in $y,y'$ and $q_y$ is taken into account by \hyp{3}. Without this argument, we would face a difference of $F$-terms which is not a priori small. A more detailed argument is given in the appendix.

We end the proof by remarking that we have a Strong Comparison Result for the $\Fb$ equation since $\Fb$ is continuous and just depends on $(p_x,p_y)$ and this allows to compare the subsolution $\oU$ and the supersolution $\uU$ to get the complete answer (see, for example, Ley\cite{ley}). This point may be a difficulty in more general cases and we refer the reader to Section~\ref{ER} for remarks in this direction.

\section{Homogenization of the $\ue/\e$ - equation}\label{Homue}

In this section, we consider (\ref{HJe}) and even a more general time dependent one, namely
\begin{equation}\label{HJet}
\ue_t + H(\e^{-1}x, \e^{-1}\ue, \e^{-1}t , D\ue)=0\quad \hbox{in}\  \R^n \times (0,+\infty),
\end{equation}
together with the initial data
\begin{equation}
\label{idHJ}
\ue (x,0)=u_0(x)\quad \hbox{in}\  \R^n .
\end{equation}
where $u_0$ is bounded, uniformly continuous in $\R^n$ and the Hamiltonians $H \in C ( \R^n \times \R \times \R^n)$ satisfies the assumptions

\medskip
\noindent \hyp{4} $H (x,u,t,p)$ is a locally Lipschitz function in $\R^n \times \R \times \R \times \R^n$, $\Z^n$-periodic in $x$, $1$ periodic in $u$ and $t$. Moreover there exists a constant $C$ such that, for almost every $(x,u,t,p)$ in $\R^n \times \R \times \R \times \R^n$
$$ |D_t H (x,u,t,p)|\leq C (1+ |H(x,u,t,p)|)\;,\;|D_u H (x,u,t,p)|\;,\; |D_p H (x,u,t,p)|\leq C.$$

\medskip
\noindent\hyp{5} $H$ is coercive in $p$, i.e.
$$H(x,u,t,p)\to +\infty\quad \hbox{as}\ |p|\to +\infty\quad \hbox{uniformly w.r.t. }\ x \in \R^n,\; u,t\in \R .$$

\medskip In addition to these assumptions, we have to add a ``geometrical'' assumption which ensures that (\ref{HJet}) is the graph-equation associated to a level-set equation in $\R^{n+1}$. A priori this requires only

\noindent\hyp{6w} There exists an Hamiltonian $H_\infty$ satisfying \hyp{4}-\hyp{5} such that
$$ H_\infty (x,u,t,p) = \lim_{s\downarrow 0}\, s H (x,u,t,s^{-1} p) \; .$$
but we use below a slightly more restrictive one

\noindent\hyp{6s} There exists a constant $C$ such that  $|D_p H (x,u,t,p)\cdot p - H (x,u,t,p)|\leq C$ a.e. in  $\R \times \R^n \times \R^n.$ \\
``w'' is for weak and ``s'' is for strong. Below we just use the notation \hyp{6} for \hyp{6s} but keeping in mind the notation $H_\infty$ introduced in \hyp{6w}.

Before providing the result, we point out that both \hyp{6w} and \hyp{6s} hold for $H(x,t,u,p) = c(x,t)|p| + g(u,t)$ if $c$, $g$ are Lipschitz, periodic functions with $c(x,t) \geq \eta >0$ in $\R^n$. Notice that $g$ may change sign or be negative.

The result is the following
\begin{theorem}
 \label{HomRes}
Assume that $H$ satisfies \hyp{4}-\hyp{6}. Then the sequence $(\ue)_\e$ converges locally uniformly to a function $u$ which is the unique solution of an equation of the type
\begin{equation}\label{HJhom}
u_t + \bar H(Du )=0\quad \hbox{in}\  \R^n \times (0,+\infty),
\end{equation}
with 
\begin{equation}
\label{idHJhom}
u (x,0)=u_0(x)\quad \hbox{in}\  \R^n .
\end{equation}
The Hamiltonian $\bar H$ can be interpreted as an homogenized Hamiltonian for a standard homogenization problem in dimension $n+1$ and is characterized by the following property: $\lambda = - \bar H (p)$ is the unique constant such that, if $w$ is the solution of 
\begin{eqnarray}
w_t + H(x,w,t,Dw)=0\quad \hbox{in}\  \R^n \times (0,+\infty),\label{eqnw}\\
w (x,0) = p\cdot x + w_0(x) \quad \hbox{in}\  \R^n ,\label{idw}
\end{eqnarray}
where $w_0$ is a bounded uniformly continuous function, then $w (x,t) - p\cdot x - \lambda t$ is uniformly bounded in $\R^n \times (0,+\infty)$.
In particular
$$ - \bar H (p) = \lim_{t\to + \infty}\, t^{-1} w(x,t)\; ,$$
the limit being uniform in $\R^n$.
\end{theorem}
  
\bigskip \noindent{\bf Proof :}\\
\noindent{\bf 1.} From $\R^n$ to $\R^{n+1}$: we interpret the equation for $\ue$ as the equation of the motion of a graph and we introduce the level-sets equation in $\R^{n+1}$. To do so, we use (a slight modification of) the approach used by Giga \& Sato \cite{GS} for Hamilton-Jacobi Equations (see also Biton, Ley and the author \cite{BBL} for the Mean Curvature Equation). The function $\Ue (x,y,t):=\ue (x,t)-y$ solves
$$ \Ue_t + |D_y \Ue| H(\e^{-1}x, \e^{-1}\Ue + \e^{-1}y, \e^{-1}t, |D_y \Ue| ^{-1}D_x \Ue)=0\quad \hbox{in}\  \R^{n+1} \times (0,+\infty),$$
$$\Ue (x,y,0)=u_0(x) - y \quad \hbox{in}\  \R^{n+1} .$$

In the sequel, it is convenient to use the notation $X=(x,y)$, $P = (p_x,p_y)$ and to introduce the Hamiltonian $F$ defined by
$$ F(X,t,P):= 
\begin{cases}
|p_y| H( x, y, t, |p_y|^{-1}p_x) & \hbox{if $p_y\neq 0$,}\\
H_\infty (x, y, t, p_x) & \hbox{otherwise.}
\end{cases}$$

By the assumptions \hyp{4}-\hyp{6} on $H$, tedious but straightforward computations show that $F$ satisfies \hyp{1}-\hyp{3} (\hyp{6s} provides, in particular, the Lipschitz continuity in $p_y$). It is worth pointing out that, in \hyp{2}, since $F$ is homogeneous of degree $1$, the coercivity of $F$ in $p_x$ or the stronger property required in the $t$-dependent case are, in fact, equivalent.

We finally consider $\Ce (X,t):= \1_{\{\Ue (\cdot,t)\geq 0\}}$ where, here and below, $\1_A$ denotes the indicator function of the set $A$. Following \cite{BSS}, it is easy to see that it solves
$$ \Ce_t + F(\e^{-1}X, \e^{-1}t, D\Ce)=0\quad \hbox{in}\  \R^{n+1} \times (0,+\infty),$$
$$\Ce (x,0)=\1_{\{u_0(x) - y\geq 0\}} \quad \hbox{in}\  \R^{n+1} .$$
Indeed the result is obvious if either $\Ue (x,t) >0$ or $\Ue (x,t)<0$ because $\Ce$ is constant in a neighborhood of $(x,t)$ and the equation is a trivial consequence of the fact that $F(X,t,0)= 0$ for any $X$ and $t$. Hence we are left with the case when $\Ue (x,t) = 0$, where the arguments of \cite{BSS} provide the answer by using the continuity of $\Ue$. In fact, in the equation, the dependence in $\Ue$ is treated as a dependence in $X$ and $t$ ; the change of variables of \cite{BSS} is done only at the level of the derivatives. This explains why the dependence in $\Ue$ disappears from this equation.

\noindent{\bf 2.} Homogenization in $\R^{n+1}$. The $\Ce$-equation can be homogenized by using almost readily Section~\ref{Hom}, the only slight difficulty being the discontinuous initial data and solutions. In fact, if $\oC:= \limsup^* \Ce$, $\uC:= \liminf_* \Ce$, $\ou:= \limsup^* \ue$, $\uu:= \liminf_* \ue$, we have the following: $\oC$ and $\uC$ are respectively viscosity sub and supersolutions of 
\begin{equation}\label{LSA}
 \chi_t + \Fb( D\chi )=0\quad \hbox{in  }\R^{n+1} \times (0,+\infty),
\end{equation}
with
$$\oC (X,0) \leq \1_{\{u_0(x) - y\geq 0\}}\; , \; \uC (X,0) \geq \1_{\{u_0(x) - y > 0\}} \quad \hbox{in}\  \R^{n+1} ,$$
and
$$\oC (X,t) = \1_{\{\ou(x,t) - y\geq 0\}}\; , \; \uC (X,t) = \1_{\{\uu(x,t) - y > 0\}} \quad \hbox{in}\  \R^{n+1}\times (0,+\infty) .$$
We also point out that $\Fb$ is homogeneous of degree $1$.\\

\noindent{\bf 3.} Reconstruction of the result in $\R^n$. Let $u: \R^{n}\times (0,+\infty) \to \R$ be the solution of 
$$ u_t + \Fb (D_x u, -1) = 0 \quad \hbox{in}\  \R^{n} \times (0,+\infty),$$
$$u (x,0)= u_0(x)  \quad \hbox{in}\  \R^{n} .$$
By the arguments of Soner, Souganidis and the author \cite{BSS}, since $u(x,t) - y$ is the solution of the level set equation (\ref{LSA}) associated to the initial data $u_0(x)-y$, we have
$$ \uC (X,t) \geq \1_{\{u(x,t) - y > 0\}} \; \hbox{   and   }\;  \oC (X,t) \leq \1_{\{u(x,t) - y \geq 0\}} \quad \hbox{in}\  \R^{n+1}\times (0,+\infty) .$$
Indeed, it is proved in \cite{BSS} that $\1_{\{u(x,t) - y > 0\}}$ and $\1_{\{u(x,t) - y \geq 0\}}$ are respectively the minimal supersolution and the maximal subsolution associated to the initial data $\1_{\{u_0(x) - y\geq 0\}}$.
By the form of $\uC$ and $\oC$ given in Step~2, this yields
$$ u(x,t) \leq \uu (x,t) \leq \ou (x,t) \leq u(x,t)  \quad \hbox{in}\  \R^{n}\times (0,+\infty) .$$
Therefore $\uu=\ou=u$ in $\R^{n}\times (0,+\infty)$ and this proves the first part of the result with $\bar H (p):= \Fb (p,-1)$.

\noindent{\bf 4.} To prove the last part of the result, we apply the same approach : if $ U : \R^{n+1}\times (0,+\infty) \to \R$ is the solution of the level set equation 
\begin{equation}\label{LSAbis}
U_t +F(X , t , DU)=0\quad \hbox{in  }\R^{n+1} \times (0,+\infty),
\end{equation}
with initial data $p\cdot x + w_0 (x) -y$, we first have $\{(x,y,t);\,w(x,t) = y\} \subset \{(x,y,t);\,U(x,y,t) = 0\}$.

On the other hand, if $v_1$ and $v_2$ are bounded space-time periodic functions which are respectively an usc subsolution and a lsc supersolution of (\ref{CP}) associated to $P=(p_x,-1)$, then the functions $v_1(X,t) + P\cdot X - \Fb (P)t$, $v_2(X,t) + P\cdot X - \Fb (P)t$ are respectively sub and supersolution of (\ref{LSAbis}); since, we may assume in addition that  $\max_{\R^{n+2}}v_1(X,t) \leq w_0 \leq \min_{\R^{n+2}} v_2(X,t)$ in $\R^{n+1}$, we have by using a comparison result
$$ v_1 (X,t) + P\cdot X - \Fb (P)t \leq 
U(X,t) \leq 
v_2 (X,t) + P\cdot X - \Fb (P)t \quad\hbox{in  }\R^{n+1} \times (0, +\infty)\; .$$
Notice that, since $H$ satisfies \hyp{4} and \hyp{6}, $F$ is Lipschitz continuous in $P$, uniformly in $X$ and we can use a property of ``finite speed of propagation'' for the level set equation, which simplifies the comparison arguments and avoids problem with the unboundedness of the domain.

If $w(x,t) = y$, then $U(x,t) = 0$, $P\cdot X = p\cdot x -w(x,t)$ and this inequality gives
$$  v_1 (X,t) + p\cdot x -w(x,t) - \Fb (P)t \leq 
0 \leq 
v_2 (X,t) + p\cdot x -w(x,t) - \Fb (P)t \quad\hbox{in  }\R^{n} \times (0, +\infty)\; ,$$
or equivalently
$$  v_1 (X,t)  \leq 
w(x,t) -p\cdot x + \Fb (P)t \leq 
v_2 (X,t) \quad\hbox{in  }\R^{n} \times (0, +\infty)\; .$$
Since $\bar H (p)= \Fb (p,-1)$ and $v_1$, $v_2$ are bounded, this gives the first part of the result. Then, dividing by $t$ and letting $t\to \infty$, we obtain the second one. And the proof is complete.

\section{Extensions and Remarks}\label{ER}

The results of the preceding sections can be extended in several directions.

First, in Section~\ref{Hom}, we can consider \textit{stationary homogenization problems}. In fact, as we suggest it in the introduction, these problems are simpler provided that we have a suitable dependence in $\ue$.

Next, the one-dimensional variable $y$ can be replaced by a \textit{multi-dimensional variable $(y_1, \cdots ,y_m)$} provided that one has the right structure on $G$ (and therefore on $F$), namely for $G$

\noindent\hyp{3'} The function $G$ is locally Lipschitz continuous w.r.t. $y=(y_1, \cdots ,y_m)$, $t$, $p_y$ and there exists $l=(l_1, \cdots, l_m) \in \R^m$ and $C_3, C_4, C_5\geq 0$ such that, for all $i=1,\cdots, m$ and for almost every $\xi = (x,y,t,p_x,p_y) \in \R^{2(n+m)+1}$
$$  |D_{y_i} G (\xi) | \leq C_3 |p_{y_i} + l_i |\; ,\; |D_t G (\xi) | \leq C_4 (1+ |p_y| + |G (\xi)|)\; , \;  |D_{p_y} G (\xi) | \leq C_5 .$$

\smallskip
\noindent With this assumption, the proof of Theorem~\ref{thm:EP} remains the same : the only slight difference is that we have to do the estimates for all the $D_{y_i}(\wa (x,y,t)+l_i y_i)$ and the first inequality for $|D_{y_i} G (\xi) |$ above implies that we can do it separately for each $i$.

For the homogenization problems, we may have \textit{additional dependences in $x,y,t$} and not only in the fast variables $\e^{-1}x, \e^{-1}y, \e^{-1}t$. In this case, $x,y,t$ become parameters in the ergodic-cell problem as are $p_x, p_y$. The main additional difficulties consists, on one hand, in obtaining the right dependence of $\Fb (x,y,t,p_x, p_y)$ in $x,y,t$ in order to have a Strong Comparison Result for the limiting equation and, on the other hand, in the convergence proof, to take care of the lack of uniform continuity of $F$ in $x, y, t$. It is worth pointing out that this last difficulty arises since we do not have (a priori) Lipschitz continuous solutions for the ergodic-cell problem. 

For the first difficulty, since one can prove, in an analogous way as in the proof of Lemma~\ref{lem:Fb}, that $\Fb (x,y,t,p_x, p_y)$ is continuous, then the dependence in $x$ does not create any problem. But the dependence in $t$ and $y$ seems less clear and we do not see how to show that (typically) $| D_y \Fb (x,y,t,p_x, p_y)| \leq C (1 +|p_y|)$ and an analogous property for the $t$-dependence. Here the fact not to have a Lipschitz continuous solution of the ergodic-cell problem creates a difficulty and we do not know how to solve it.

On the contrary, the second difficulty can be solved by using a similar trick to the $F^k$-one and is not a real problem.

Finally, even if, in this paper, we restrict ourselves to the case of first-order equations, there are easy extensions to second-order equations~: for example, Theorem~\ref{thm:EP} extends without any difficulty if we add a $u_{yy}$ to the equation, which corresponds to homogenization problems with a vanishing viscosity term of the form $\e^2\ue_{yy}$. Such type of extensions will be considered in a future work.

\section*{Appendix : Proof of (\protect\ref{fkvi})}

Since $(\xb, \yb, \tb)$ is a {\sl strict} maximum point of $\oU - \phi$, standard arguments show that, for $\e$ and $\delta$ small enough, the function
$$ \Ue (x,y,t) - \phi (x,y,t) - \e V (\e^{-1}x', \e^{-1}y',\e^{-1}t') - \frac{|x-x'|^2}{\delta^2} - \frac{|y-y'|^2}{\delta^2}- \frac{|t-t'|^2}{\delta^2}\; ,$$
has a local maximum point near $(\xb, \yb, \tb,\xb, \yb, \tb)$. We denote it by $(x,y,t,x',y',t')$ for the sake of notational simplicity. Moreover we know that $(x,y,t,x',y',t')$ converges to $(\xb, \yb, \tb,\xb, \yb, \tb)$ as $\e, \delta \to 0$ and
$$ \frac{|x-x'|^2}{\delta^2} + \frac{|y-y'|^2}{\delta^2} + \frac{|t-t'|^2}{\delta^2} \to 0\quad\hbox{as  $\delta \to 0$ for any $\e >0$.}$$

We set
$$ q_x := \frac{2(x-x')}{\delta^2}\; ,  \; q_y := \frac{2(y-y')}{\delta^2}\; ,  \; q_t := \frac{2(t-t')}{\delta^2}\; .$$
The viscosity inequalities for $\Ue$ at $(x,y,t)$ and for $V$ at $(\e^{-1}x', \e^{-1}y',\e^{-1}t')$ (after applying the obvious unscaling) read
$$ \phi_t (x,y,t) + q_t + F(\e^{-1}x, \e^{-1}y, \e^{-1}t, D_x \phi (x,y,t) + q_x, D_y  \phi (x,y,t)+ q_y) \leq 0\; ,
$$ 
$$
q_t + F^k (\e^{-1}x', \e^{-1}y', \e^{-1}t' ,q_x + p_x , q_y + p_y ) \geq \Fb^k (p_x ,p_y)
$$
where we recall that $p_x = D_x \phi (\xb, \yb, \tb)$ and $p_y = D_y \phi (\xb, \yb, \tb)$.

But, by the above recalled property on the penalization terms, we have
$$ |\e^{-1}x-\e^{-1}x'|^2 + |\e^{-1} t-\e^{-1} t'| = \e^{-1}\delta o(1)\quad \hbox{where $o(1) \to 0$ as  }\e, \delta \to 0\; ,$$
and, on the other hand, since $(x,y,t) \to (\xb, \yb, \tb)$, we have
$$ |(D_x \phi (x,y,t) + q_x) - (q_x + p_x)| \to 0\; \hbox{as  }\e, \delta \to 0\; .$$
Therefore if we choose $\delta \ll \e$ and $\e, \delta$ small enough, we have
\begin{equation}\label{inxx}
|\e^{-1}x-\e^{-1}x'|^2 , |\e^{-1} t-\e^{-1} t'|, |(D_x \phi (x,y,t) + q_x) - (q_x + p_x)| \leq k^{-1}\; .
\end{equation}

Next, we use this information together with \hyp{3} to obtain
$$
F(\e^{-1}x, \e^{-1}y, \e^{-1}t, D_x \phi (x,y,t) + q_x, D_y  \phi (x,y,t)+ q_y)  \geq $$
$$ F(\e^{-1}x, \e^{-1}y', \e^{-1}t, D_x \phi (x,y,t) + q_x, D_y  \phi (\xb, \yb, \tb) + q_y)  - C_3 |\e^{-1}y-\e^{-1}y'| |D_y  \phi (x,y,t)+ q_y| $$ $$- C_5 |D_y  \phi (x,y,t)-D_y  \phi (\xb, \yb, \tb)| \; .$$
And thanks again to the property on the penalization terms, we can choose $\delta$ small enough in order to have
$$ |\e^{-1}y-\e^{-1}y'| |D_y  \phi (x,y,t)+ q_y| = o(1) \quad\hbox{as  }\e, \delta \to 0\; .$$
Using this and (\ref{inxx}) yields the final estimate of $F(\e^{-1}x, \e^{-1}y, \e^{-1}t, D_x \phi (x,y,t) + q_x, D_y  \phi (x,y,t)+ q_y)$, namely
$$ F^k(\e^{-1}x', \e^{-1}y', \e^{-1}t, D_x \phi (\xb, \yb, \tb) + q_x, D_y  \phi (\xb, \yb, \tb) + q_y)  - o(1)\; .$$

In order to conclude, it suffices to subtract the above viscosity inequalities and to use this information~: this yields
$$ \phi_t (x,y,t) -o(1) \leq  - \Fb^k (p_x ,p_y) \; ,
$$
and we let $\e, \delta$ tends to $0$, choosing $\delta$ sufficiently small compared to $\e$.

\begin{Remark} This detailed proof shows why the ``$F^k$-trick'' can be useful : without using it, we would have to compare the terms $F(\e^{-1}x, \e^{-1}y, \e^{-1}t, D_x \phi (x,y,t) + q_x, D_y  \phi (x,y,t)+ q_y)$ and $F (\e^{-1}x', \e^{-1}y', \e^{-1}t' ,q_x + p_x , q_y + p_y )$. As the above proof shows it, the ``good'' assumption \hyp{3} allows to treat the differences in $y$ and $D_y u$ and it seems to be a natural assumption to do it. But, for the $x$ and $t$ variables, we do not have such type of assumption and the ``$F^k$-trick'' seems the only way to handle the corresponding differences. 

We have choosen to present this ``$F^k$-trick'' in this way in order to point out the difference between having \hyp{3} and not having it, but we could have defined $F^k$ in a different way in order that it takes care of all the variables.
\end{Remark}

\end{document}